\documentclass[leqno,draft]{article}

\def\diam{\mathop{\rm diam}}
\def\dist{\mathop{\rm dist}}

\begin{document}

\title{Complementary self-similarity}

\author{Stephen Semmes \\
        Rice University}

\date{}

\maketitle

        Let $U$ be a nonempty proper open subset of ${\bf R}^n$, and
put $E = \partial U$.  Thus $E$ is a nonempty closed set in ${\bf
R}^n$ with empty interior, and we are especially interested in the
case where $E$ is connected and $U$ has infinitely many connected
components.  For instance, it may be that $U = {\bf R}^n \backslash
E$, or that $U$ is the union of the bounded components of the
complement of $E$.

        A strong version of ``complementary self-similarity'' of $E$
would ask that the components of $U$ be similar to each other, which
is to say that they can be transformed into one another by
combinations of translations, orthogonal transformations, and
dilations.  Many standard examples such as Sierpinski gaskets and
carpets have this property, since the bounded complementary components
are equilateral triangles or squares, etc.  There are obviously
numerous extensions of this, using bilipschitz, quasiconformal, or
quasisymmetric mappings, and so on.  However, this does not
necessarily include self-similarity of the arrangement of the
components of $U$.  As an extreme case, it may be that $E$ has
positive Lebesge measure, and that the sum of the $(n -
1)$-dimensional Hausdorff measures of the boundaries of the components
of $U$ is finite.

        One could also consider scale-invariant geometric conditions
on the components of $U$.  A very basic condition of this type would
ask that each component $V$ of $U$ contain an open ball whose radius
is at least a positive constant times the diameter of $V$.  A related
scale-invariant condition concerning the arrangement of the components
of $U$ would say that for each $x \in E$ and $0 < r \le \diam U$ there
is a component of $U$ contained in $B(x, r)$ whose diameter is at
least a positive constant times $r$.  The combination of these two
conditions implies that $E$ is \emph{porous}, in the sense that for
each $x \in E$ and $0 < r \le \diam U$ there be a point in $B(x, r)
\backslash E$ whose distance to $E$ is at least a positive constant
times $r$.  A more qualitative form of the second condition would
simply ask that there be a component of $U$ contained in $B(x, r)$ for
every $x \in E$ and $r > 0$.

        As another possibility, suppose that there is a $k \ge 1$ such
that for each component $V$ of $U$ and every $x, y \in \partial V$
there is a continuous path in $\partial V$ that connects $x$ to $y$
and has length $\le k \, |x - y|$.  This implies that every continuous
path in $\overline{U}$ with finite length $\le L$ connecting a pair of
elements of $E$ can be modified to get a path in $E$ with the same
endpoints and length $\le k \, L$.  In many cases, $\overline{U}$ is
very nice, perhaps even convex, and so it may be very easy to connect
elements of $\overline{U}$ by paths of controlled length.  Of course,
there are analogous assertions for higher-dimensional fillings.

        Although distinct components of $U$ are disjoint, their
closures may or may not be.  If the closures are disjoint, then one
can look for lower bounds on the distances between them.  A basic
scale-invariant condition would ask that
\begin{equation}
        \min(\diam V, \diam W) \le C \, \dist(V, W)
\end{equation}
for some $C > 0$ and any two distinct components $V$, $W$ of $U$.

        Note that $E$ has topological dimension $n - 1$ when $E =
\partial U$ for a disconnected open set $U \subseteq {\bf R}^n$.  More
precisely, the topological dimension of $E$ is $\le n - 1$ because $E$
has empty interior, and it is then equal to $n - 1$ because $E$
separates elements of ${\bf R}^n$.

        Suppose that $E$ is compact.  For each $w \in {\bf R}^n
\backslash E$, let $\pi_w$ be the mapping from $E$ into the unit
sphere ${\bf S}^{n - 1}$ defined by
\begin{equation}
        \pi_w(x) = \frac{x - w}{|x - w|}.
\end{equation}
A famous theorem states that $\pi_w$ and $\pi_z$ are homotopically
equivalent if and only if $w$ and $z$ are in the same component of the
complement of $E$.  In particular, $\pi_w$ is homotopically trivial
exactly when $w$ is in the unbounded component of ${\bf R}^n
\backslash E$.  The Lipschitz constant of $\pi_w$ can be estimated in
terms of $1/\dist(w, E)$, which one might like to be as small as
possible.

        Even when the complement of $E$ is connected, there can be a
lot of topological activity related to linking.  There is also some
discreteness involved with this, as indicated by homology or
cohomology groups with integer coefficients, or by homotopy classes
into nice spaces for which close mappings are homotopic.  There can
still be matters of linking when $E$ is connected, $U$ has infinitely
many connected components, and $n \ge 3$, because the components of
$U$ may not be contractable.  Common features of the topology of the
components of $U$ are another aspect of self-similarity of the
complement of $E$.  At any rate, one is likely to be interested in
more than just homotopy classes of continuous mappings from $E$ into
the nonzero complex numbers when the topological dimension of $E$ is
strictly greater than $1$.

        A nice hypersurface in ${\bf R}^n$ has two sides,
corresponding to two complementary components.  Using Cauchy integrals
from ordinary complex analysis when $n = 2$ and extensions to higher
dimensions based on quaternions or Clifford algebras, functions on the
hypersurface can be split into two pieces, which are boundary values
of holomorphic functions on the two complementary components.
However, if $E \subseteq {\bf R}^n$ is sufficiently ``large'', then
the Cauchy integral of suitable measures on $E$ determine continuous
functions on ${\bf R}^n$ that are holomorphic on the complement of
$E$.  One can still look at operators involving projections onto
holomorphic functions on the individual complementary components, even
when there are infinitely many components.

\end{document}